\documentclass[11pt, final]{amsart}
\usepackage{amssymb}
\usepackage[mathscr]{eucal}

\theoremstyle{plain}
\newtheorem{theorem}{Theorem}[section]
\newtheorem{proposition}[theorem]{Proposition}
\newtheorem{lemma}[theorem]{Lemma}

\theoremstyle{remark}
\newtheorem*{remark}{Remark}

\def\cA{\mathscr{A}}
\def\cB{\mathscr{B}(\mathscr{H})}
\def\cH{\mathscr{H}}
\def\C{\mathbb C}
\def\CX{C(X, \cB)}
\def\d{\, d}
\def\Jaut{\operatorname{Jaut}}

\def\M{\mathscr{M}}
\def\N{\mathbb N}
\def\refl{\operatorname{ref}}
\def\R{\mathbb R}

\begin{document}
\title[]{A reflexivity problem concerning
the $C^*$-algebra $C(X)\otimes \cB$}
\author{LAJOS MOLN\'AR}
\address{Institute of Mathematics\\
         Lajos Kossuth University\\
         4010 Debrecen, P.O.Box 12, Hungary}
\email{molnarl@math.klte.hu}
\thanks{  This research was supported from the following sources:\\
          1) Joint Hungarian-Slovene research project supported
          by OMFB in Hungary and the Ministry of Science and
          Technology in Slovenia, Reg. No. SLO-2/96,\\
          2) Hungarian National Foundation for Scientific Research
          (OTKA), Grant No. T--030082 F--019322,\\
          3) A grant from the Ministry of Education, Hungary, Reg.
          No. FKFP 0304/1997}
\begin{abstract}
Let $X$ be a compact Hausdorff space and let $\cH$ be a separable
Hilbert space. We prove that the group of all order automorphisms of the
$C^*$-algebra $C(X)\otimes \cB$ is algebraically reflexive.
\end{abstract}
\maketitle

\section{Introduction}
Reflexivity problems concerning subspaces of the algebra $\cB$ of all
bounded linear
operators acting on a Hilbert space $\cH$ is one of the most active
research areas in operator theory. In the past decade,
similar questions concerning some important subsets of
linear transformations acting on Banach algebras (not on Hilbert spaces)
have also attracted attention.
The originators of the research in this direction are Kadison and
Larson. In the papers \cite{Kad}, \cite{LarSour} the reflexivity of the
Lie algebra of
derivations was treated. In \cite[Some concluding remarks (5),
p. 298]{Lar} Larson raised the reflexivity problem for automorphism
groups. This was investigated together with the similar problem for
isometry groups in a series of papers
\cite{BaMo, MolStud1, MolGyor, MolZal1, MolZal2, MolLond}.
For example, in \cite{MolGyor} we proved that the automorphism group and
the isometry group of the suspension of $\cB$, $\cH$ being a separable
infinite dimensional Hilbert space, are algebraically reflexive.
As for the automorphism group, this result was a consequence
of \cite[Theorem 2]{MolGyor} stating that for any locally compact
Hausdorff space $X$, if the automorphism group of $C_0(X)$ (the algebra
of all continuous complex valued functions on $X$ vanishing at infinity)
is reflexive, then so is the automorphism group of $C_0(X)\otimes \cB$.
The main theorem of the first part of the present paper gives a similar
result for the group of all order automorphisms in case $X$ is
compact. The $C^*$-algebra $C(X)\otimes \cB$
is of importance because of several remarkable
results concerning the structure of its derivations and automorphisms
\cite{Lan} (see also \cite{KadRing}).

Let us fix the concepts and the notation what we use throughout.
Let $\mathscr X$ be a Banach space (in fact, we shall be mostly
interested in the case when it is a $C^*$-algebra) and for any subset
$\mathscr E\subset {\mathscr B}(\mathscr X)$ we define
\[
\refl_{al} \mathscr E=\{ T\in {\mathscr B}(\mathscr X) : Tx\in \mathscr
Ex \text{ for all } x\in \mathscr X\}.
\]
The elements of $\refl_{al} \mathscr E$ can be described as those linear
transformations which are locally in  $\mathscr E$. The set $\mathscr E$
is called algebraically reflexive if $\refl_{al}\mathscr E=\mathscr E$,
that is, when, so to say, every linear transformation
locally in $\mathscr E$ is globally in $\mathscr E$.
We use the expression ``algebraically reflexive'' since there is an
analogue concept called ``topological reflexivity''.

If $\cA$ is a $C^*$-algebra, denote by $\cA^+$ the cone of all
positive elements in $\cA$. The linear bijection $\phi:\cA \to
\cA$ is called an order automorphism of $\cA$ if $\phi(a)$ is positive
if and only if $a$ is positive $(a\in \cA)$. Clearly, every linear map
$\psi$ on $\cA$ which is a local order automorphism (meaning that $\psi$ is
locally
in the set of all order automorphisms of $\cA$) preserves the positivity
in both directions ($\psi(a)\geq 0$ if and only if $a\geq 0$) and it
is injective. Therefore, to prove that the goup of order automorphisms
of $\cA$ is algebraically reflexive, the only thing which needs
verification is that every local order automorphism of $\cA$ is
surjective.

If $\mathscr{R}$ and $\mathscr{R'}$ are *-algebras, then the linear map
$\phi:\mathscr R \to \mathscr R'$ is called a Jordan *-homomorphism if
it satisfies
\[
\phi(x^2)=\phi(x)^2, \qquad \phi(x^*)=\phi(x)^* \quad (x \in \mathscr
R).
\]
The definition of Jordan *-isomorphisms and that of Jordan
*-automorphisms should be clear. The set of all Jordan *-automorphisms
of the *-algebra $\mathscr R$ is denoted by $\Jaut(\mathscr R)$.

Let $X$ be a compact Hausdorff space and let $\cH$ be a Hilbert space.
Consider the tensor product $C(X)\otimes \cB$ which is meant in the
$C^*$-algebra
sense. It is well-known that this $C^*$-algebra is isomorphic and
isometric to the algebra $\CX$ of all continuous functions from $X$
into $\cB$ and that the linear span of the elements $fA$ $(f\in C(X),
A\in \cB)$ is dense in $\CX$. The unit of $C(X)$ is denoted by $1$ while
the unit of $\cB$ is denoted by $I$.

\section{The Results}

Our first theorem describes the order automorphisms of $\CX$.

\begin{theorem}\label{T:formord}
Let $\cH$ be a separable Hilbert space and let $X$
be a compact Hausdorff space.
If $\phi$ is an order automorphism of $C(X, \cB)$, then $\phi$ is
of the form
\begin{equation}\label{E:formpos}
\phi(f)(x)=b(x)([\tau(x)](f(\varphi(x))))b(x)^*
\qquad (f\in C(X, \cB), x \in X),
\end{equation}
where $b\in C(X, \cB)$ is such that its values are invertible
positive operators, $\tau: X \to \Jaut(\cB)$ is a function such that
$x \mapsto \tau(x)$ and $x \mapsto \tau(x)^{-1}$ are strongly continuous
and $\varphi :X \to X$ is a homeomorphism.
\end{theorem}

\begin{proof}
First suppose that $\phi$ is unital, that is, $\phi(1I)=1I$.
According to a well-known theorem of Kadison
\cite[Corollary 5]{Kad2} every unital order automorphism of a
$C^*$-algebra is a Jordan *-automorphism. Thus $\phi$ is a
Jordan *-automorphism of
$C(X, \cB)$. Since the closed Jordan ideals of a $C^*$-algebra coincide
with its closed associative ideals, following the arguments
used in the proofs of \cite[Lemma 2.1 and Theorem 1]{MolGyor}, we obtain
the form \eqref{E:formpos} with $b=1I$. Otherwise, define
$c(x) =\sqrt{\phi(1I)(x)}^{-1}$. We show that the unital map $\psi: C(X,
\cB)\to C(X, \cB)$ defined by
\[
\psi(f)(x)=c(x)\phi(f)(x)c(x) \qquad
(f\in C(X, \cB), x \in X)
\]
is an order automorphism. To see this, we only have to prove that
$c:X \to \cB$ is continuous. But this follows from the continuity of the
inverse and the square-root operations.
It is now trivial to complete the proof.
\end{proof}

\begin{remark}
In fact, \eqref{E:formpos} can be written in a more
explicit form if one reminds that the Jordan *-automorphisms of the
algebra $\cB$ are exactly the maps
\[
A \longmapsto UAU^*, \qquad A \longmapsto VA^{tr} V^*,
\]
where $U, V$ are unitary operators on $\cH$ and ${}^{tr}$ denotes the
transpose with respect to an arbitrary but fixed orthonormal basis in
$\cH$.
\end{remark}

Concerning the reflexivity of the group of all order automorphisms of
the function algebra $C(X)$ we have the following result.
Using the above mentioned result of Kadison, for example,
one can easily verify that the order automorphisms of $C(X)$ are
precisely the maps
\[
f \longmapsto h \cdot f \circ \varphi,
\]
where $h\in C(X)$ is everywhere positive and $\varphi :X \to X$ is a
homeomorphism.

\begin{theorem}\label{T:ordrefcx}
Let $X$ be a first countable compact Hausdorff space.
Then the group of all order automorphisms of $C(X)$ is algebraically
reflexive.
\end{theorem}

\begin{proof}
Let $\phi :C(X) \to C(X)$ be a linear map which is a local order
automorphism of $C(X)$.
Clearly, we may suppose that $\phi(1)=1$. We claim that in this case
$\phi$ is an algebra endomorphism of $C(X)$. To see this, let
$x\in X$ be fixed and consider the linear functional
$\phi_x(f)=\phi(f)(x)$ $(f\in C(X))$. Since this is a positive linear
functional, by the Riesz representation theorem we have a regular
probability measure $\mu$ on the Borel sets of $X$ such that
\[
\phi_x(f)=\int_X f \d\mu \qquad (f \in C(X)).
\]
Suppose that there are two disjoint closed sets $A,B$ in $X$ which are
of positive measure. Then, by Urysohn's lemma, we can choose continuous
functions $f,g:X \to [0,1]$ with disjoint support such that $f_{|A}=1$,
$f_{|B}=0$
and $g_{|A}=0$, $g_{|B}=1$. Consider the function $f+ig$. Clearly, its
integral is a complex number with nonzero real and imaginary parts.
On the other hand, from
the local property of $\phi$ it follows that $\phi_x(f+ig)$ should be a
positive scalar multiple of one of the values of $f+ig$. Since the
range of $f+ig$ lies in $\R \cup i \R$, we arrive at a contradiction.
This shows that either $\mu(A)=0$ or $\mu(B)=0$.
Since this holds for every pair of disjoint closed sets of $X$, by
regularity we obtain that every Borel set in $X$ has measure either 0 or
1. This implies that $\mu$ is a Dirac measure.
In fact, it follows that the integral (as a number) of every
simple function
is contained in the range of that function. Then, approximating
continuous functions by simple ones, we find that the integral of
every continuous function $f \in C(X)$ belongs to its range. Therefore,
$\phi_x$ is a linear functional on $C(X)$ which sends 1 to 1 and has the
property that $\phi_x(f)$ belongs to the spectrum of $f$. By the famous
Gleason-Kahane-\. Zelazko theorem it follows that $\phi_x$ is
multiplicative and, therefore, a point-evaluation which gives us that
its representing measure is a Dirac measure.
Hence we obtain that $\phi$ is a unital endomorphism
of $C(X)$. But the form such morphisms of $C(X)$ is
well-known. Namely, there is a continuous function $\varphi: X\to X$
such that
\[
\phi(f)=f \circ \varphi \qquad (f\in C(X)).
\]
It remains to prove that $\varphi$ is bijective. The surjectivity of
$\varphi$ follows from the injectivity of $\phi$. To the injectivity of
$\varphi$ observe that using Urysohn's lemma, by the first countability
of $X$, we have a nonnegative function $f_0$ on $X$ which vanishes at
exactly one point. If $\varphi$ was noninjective, we would obtain that
$\phi(f_0)$ vanishes at more than one point. But this is a contradiction
due to the local property of $\phi$. This completes the proof.
\end{proof}

\begin{remark}
We show that in our previous result the condition of first countability
is essential. The following example was communicated to us by
F\'elix Cabello S\'anchez on the occasion of another reflexivity
problem. Consider the compact Hausdorff space $\N^*=\beta \N \setminus
\N$, where $\beta \N$ stands for the Stone-\v Cech compactification of
the positive integers. Let $A,B \subset \N$ be disjoint infinite sets.
Then, by \cite[3.10, 3.14, 3.15]{Walker}, the sets $A^*=cl_{\beta \N}
A\setminus A$ and $B^*=cl_{\beta \N} B\setminus B$ are disjoint clopen
sets in $\N^*$ and they are both homeomorphic to $\N$. This gives us the
existence of a continuous surjective and noninjective function $\varphi:
\N^* \to \N^*$. Consider the map $\phi: f\mapsto f \circ \varphi$.
It is an interesting property of $\N^*$ that if two functions
$f,g \in C(\N^*)$ have the same range, then there is a homeomorphism
$h :\N^* \to \N^*$ such that $f=g\circ h$ \cite[first section on p.
83]{Walker}. It then follows that for every $f\in C(\N^*)$ there exists
a homeomorphism $h: \N^* \to \N^*$ such that $f\circ \varphi=f\circ h$.
Consequently, $\phi$ is a local order automorphism of $C(\N^*)$ which is
not surjective.
\end{remark}

The following result shows that there are compact Hausdorff spaces
which are not first countable but have our reflexivity property.

\begin{proposition}
The group of all order automorphisms of $C(\beta \N)\cong \ell_\infty$
is algebraically reflexive.
\end{proposition}

\begin{proof}
Let $\phi:C(\beta \N) \to C(\beta \N)$ be a linear map which is a local
order automorphism. As before, we can suppose that $\phi(1)=1$. By the
proof of Theorem~\ref{T:ordrefcx} we obtain that there is a surjective
function
$\varphi: \beta \N \to \beta \N$ such that $\phi(f)=f\circ \varphi$
$(f\in C(\beta \N))$.
We assert that $\varphi(\N) \subset \N$. Let $c_0$
denote the space of all complex sequences converging to 0. Observe that
under the identification of $\ell_\infty $ and $C(\beta \N)$, $c_0$ can
be considered as a subalgebra of $C(\beta \N)$. Since the set of all
isolated points of $\beta \N$ is exactly $\N$, every homeomorphism of
$\beta \N$ maps $\N$ onto $\N$. Let $f\in C(\beta \N)$ be defined
by $f(n)=1/n$ $(n \in \N)$. By the local property of $\phi$ it follows
that the set of all nonzeros of $\phi(f)$ is $\N$.
On the other hand, since $\phi(f)=f\circ \varphi$, we obtain that the
set of all nonzeros of $\phi(f)$ is the preimage of the set of nonzeros
of $f$ under $\varphi$. Consequently, we have $\N =\varphi^{-1}(\N)$.

Clearly, $\phi$ can be considered as a map from $\ell_\infty$ into
itself. We obtain
\[
\phi( (\lambda_n)) =(\lambda_{\varphi(n)})
\qquad ((\lambda_n)\in \ell_\infty).
\]
Here $\varphi$ is treated as a function from $\N$ into itself.
Since $\phi$ is injective, we obtain that $\varphi$ maps onto $\N$.
We are done if we prove that $\varphi$ is injective as well. To see
this, observe that for any sequence $(\lambda_n)$ with entries all zero
but one, $\phi((\lambda_n))$ must have the same property.
The proof is now complete.
\end{proof}

In the proof of our main result we use the
following easy lemma.

\begin{lemma}\label{L:possubs}
The only nontrivial subspaces in $\C \cB^+$ which contain an
invertible operator are one-dimensional.
\end{lemma}

\begin{proof}
Let $\M$ be a nontrivial subspace with the above property.
Clearly, we can assume that $I\in \M$. Let $A\in \M$ be arbitrary.
Since $A+\lambda I$ is a scalar multiple of a positive operator, it
follows that for every $\lambda \in \C$,
the spectrum $\sigma(A) +\lambda$ of $A+\lambda I$ lies on
a straight line of the plane going through 0. This trivially
implies that the spectrum of $A$ consists of one element, which, by the
normality of $A$, gives us that $A$ is a scalar.
\end{proof}

\begin{remark}
We note that the conclusion in the above lemma remains valid
even if we do not assume that the subspace in question contains
an invertible element. In that case the proof is not so trivial and can
be based on an observation of Radjavi and Rosenthal \cite[Remark
iii, p. 691]{RaRo} stating that every subspace of normal operators is
commutative. However, the above lemma is sufficient for our
purposes.
\end{remark}

\begin{theorem}\label{T:ordrefbh}
Let $\cH$ be a separable Hilbert space. If the
group of all order automorphisms of $C(X)$ is algebraically reflexive,
then so is the group of all order automorphisms of $C(X) \otimes
\cB$.
\end{theorem}

\begin{proof}
Let $\phi :\CX \to \CX$ be a linear map which is a local order
automorphism of $\CX$. Clearly, we may assume that $\phi(1I)=1I$.

Fix $x\in X$ and consider the map
\[
f\longmapsto \phi(fI)(x)
\]
on $C(X)$. By Theorem~\ref{T:formord}, this is a linear map whose values
are lying in $\C \cB^+$. Since $\phi(1I)=1I$, from Lemma~\ref{L:possubs}
it
follows that the map $f\mapsto \phi(fI)$ can be considered as a linear
transformation on $C(X)$ into itself.
We assert that it is an automorphism of $C(X)$.
By the local form of $\phi$, for every $f\in C(X)$ there exists
a function $h_f\in \CX$ whose values are invertible and positive, and
a homeomorphism $\varphi_f :X \to X$ such that
\[
\phi(fI)=f\circ \varphi_f\cdot h_f.
\]
Since the values of $\phi(fI)$ are scalars, we easily obtain that
$h_f$ can be chosen to be scalar valued as well. So, the map
\[
f \longmapsto \phi(fI)
\]
can be considered as a local order automorphism of $C(X)$. By our
assumption on $C(X)$ it follows that this map is an order automorphism
of $C(X)$. Since we have supposed that $\phi(1I)=1I$, we infer that
the above map is in fact an automorphism of $C(X)$, that is,
there exists a homeomorphism $\varphi :X \to X$ for which
\[
\phi(fI)=(f\circ \varphi) I \qquad (f \in C(X)).
\]

Let us fix a point $x\in X$ and consider the map
\[
A \longmapsto \phi(1A)(x)
\]
on $\cB$. By the local form of $\phi$, the above map is a local order
automorphism of $\cB$. By \cite[Theorem 2]{MolSem}, every local order
automorphism of $\cB$ is an order automorphism. Therefore, it follows
that there is a positive invertible operator
$b(x)\in \cB$ and a Jordan *-automorphism $\tau(x)$ of $\cB$ such that
\[
\phi(1A)(x)=b(x)[\tau(x)](A)b(x) \qquad (x\in X).
\]
Since $\phi$ is unital, we infer that $b(x)^2=I$ which gives us that
$b(x)=I$. Hence, there is a function $\tau: X \to \Jaut(\cB)$ such
that
\[
\phi(1A)(x)=[\tau(x)](A) \qquad (x\in X).
\]
Clearly, $\tau$ is strongly continuous.

We claim that
\[
\phi(fA)=f(\varphi(x))([\tau(x)](A)) \qquad
(f\in C(X), A\in \cB, x\in X).
\]
To see this, pick $x\in X$. Let $A\in \cB$ be positive and invertible.
Consider the linear map
\[
f \longmapsto \phi(fA)(x).
\]
The image of this map is a linear subspace in $\C \cB^+$ containing the
invertible operator $\phi(1A)(x)=[\tau(x)](A)$. So, by
Lemma~\ref{L:possubs} it follows that every member of the above linear
subspace is a scalar multiple of $[\tau(x)](A)$.
Fixing $f$ for a moment and considering any other positive
invertible operator $B$, we have constants $\alpha, \beta, \gamma \in
\C$ such that
\begin{gather*}
\phi(fA)(x)=\alpha [\tau(x)](A),\\
\phi(fB)(x)=\beta [\tau(x)](B),\\
\phi(f(A+B))(x)=\gamma [\tau(x)](A+B).
\end{gather*}
By the additivity of $\phi$, in case $A,B$ are linearly independent, we
find that $\alpha=\beta$. Clearly, this is the case also when $A,B$ are
linearly dependent. This yields that for every $f\in C(X)$ and $x\in
X$, there is a constant $\lambda_{f,x}$ such that
\[
\phi(fA)(x)=\lambda_{f,x}[\tau(x)](A)
\]
holds true for every positive invertible operator $A$.
Putting $A=I$ we obtain
\[
f(\varphi(x))I=\phi(fI)(x)=\lambda_{f,x}[\tau(x)](I)=
\lambda_{f,x}.
\]
This means that
\begin{equation}\label{E:1}
\phi(fA)(x)=f(\varphi(x))[\tau(x)](A)
\end{equation}
holds for every $f\in C(X)$, $x\in X$, and positive invertible
$A\in \cB$. Since every operator in $\cB$ is a linear combination of
positive invertible operators, it follows that \eqref{E:1}
holds true for every operator $A\in \cB$. This gives us that
\[
\phi(fA)(x)=[\tau(x)]((fA)(\varphi(x)))
\]
for every $f\in C(X), x\in X$ and $A\in \cB$. Since the set of all
elements of the form $fA$ $(f\in C(X), A\in \cB)$ is dense in $\CX$, by
the continuity of $\phi$ we obtain
\begin{equation}\label{E:2}
\phi(f)(x)=[\tau(x)](f(\varphi(x)))
\end{equation}
for every $f\in \CX$. Up to this point,
everything was done in order to be able to prove that $\phi$ is
surjective. This is now easy. In fact, we first verify that
$x\mapsto \tau(x)^{-1}$ is strongly continuous.
Since every Jordan *-automorphism of a $C^*$-algebra is an isometry,
this follows from the equality
\begin{gather*}
\|\tau(x)^{-1}(A)-\tau(x_0)^{-1}(A)\|=
\|A-\tau(x)(\tau(x_0)^{-1}(A))\|=\\
\|\tau(x_0)(\tau(x_0)^{-1}(A))-\tau(x)(\tau(x_0)^{-1}(A))\|
\end{gather*}
and from the strong continuity of $\tau$.
To show the surjectivity of $\phi$,
it is now enough to check that for every $g\in \CX$, the
function $x \mapsto [\tau(x)^{-1}](g(x))$ belongs to $\CX$.
But this follows from the inequality
\begin{gather*}
\| \tau(x)^{-1}(g(x))-\tau(x_0)^{-1}(g(x_0))\|=\\
\| \tau(x)^{-1}(g(x)-g(x_0))\| + \|
\tau(x)^{-1}(g(x_0))-\tau(x_0)^{-1}(g(x_0))\|\leq\\
\| g(x)-g(x_0)\| + \|
\tau(x)^{-1}(g(x_0))-\tau(x_0)^{-1}(g(x_0))\|
\end{gather*}
using the continuity of $g$ and the strong continuity of $x \mapsto
\tau(x)^{-1}$. The proof is now complete.
\end{proof}

\bibliographystyle{amsplain}

\end{document}